\documentclass[algo2e]{ifacconf}

\usepackage{graphicx}      % include this line if your document contains figures
\usepackage{natbib}        % required for bibliography
\usepackage{color}
\usepackage{amssymb}
\usepackage{amsmath}
\usepackage{algorithm}
\usepackage{algpseudocode}
\usepackage{mathtools}
\usepackage{xspace}
\usepackage{listings}
\usepackage{matlab-prettifier}
\usepackage{siunitx}
\usepackage{array}
\newcommand\Tstrut{\rule{0pt}{4mm}} 

\graphicspath{{fig/}}

\newcommand{\boldalpha}{\boldsymbol{\alpha}}
\newcommand{\boldbeta}{\boldsymbol{\beta}}
\newcommand{\boldomega}{\boldsymbol{\omega}}
\newcommand{\p}{\mathbf{p}}
\newcommand{\boldd}{\mathbf{d}}
\newcommand{\R}{\mathbf{R}}
\newcommand{\M}{\mathbf{M}}
\newcommand{\boldL}{\mathbf{L}}
\newcommand{\q}{\mathbf{q}}
\newcommand{\bolde}{\mathbf{e}}
\newcommand{\x}{\mathbf{x}}
\newcommand{\boldc}{\mathbf{c}}
\newcommand{\mbbR}{\mathbb{R}}
\newcommand{\boldu}{\mathbf{u}}

\newcommand{\calS}{\mathcal{S}}
\newcommand{\calE}{\mathcal{E}}
\newcommand{\Pf}{\mathbb{P}_f}
\newcommand{\PPP}{\mathrm{PPP}}
\newcommand{\Ell}{\mathrm{Ell}}
\newcommand{\0}{\mathbf{0}}

\newcommand{\SO}{\mathrm{SO(3)}}
\newcommand{\SE}{\mathrm{SE(3)}}
\newcommand{\Vs}{\ensuremath{\mathbb{V}^\calS}\xspace}
\newcommand{\Vt}{\ensuremath{\mathbb{V}^\mathrm{t}}\xspace}
\newcommand{\Vo}{\ensuremath{\mathbb{V}^\mathrm{o}}\xspace}
\newcommand{\Vh}{\ensuremath{\mathbb{V}^\mathrm{h}}\xspace}
\newcommand{\vs}{\ensuremath{\boldsymbol{v}^\calS}\xspace}
\newcommand{\vt}{\ensuremath{\boldsymbol{v}^\mathrm{t}}\xspace}
\newcommand{\vo}{\ensuremath{\boldsymbol{v}^\mathrm{o}}\xspace}
\newcommand{\vh}{\ensuremath{\boldsymbol{v}^\mathrm{h}}\xspace}
\newcommand{\ns}{\ensuremath{n_\mathrm{s}}\xspace}
\newcommand{\nt}{\ensuremath{n_\mathrm{t}}\xspace}

\newcommand{\noj}{\ensuremath{n­_\mathrm{oj}}\xspace}

\newcommand{\nh}{\ensuremath{n_\mathrm{h}}\xspace}
\newcommand{\cs}{\ensuremath{\boldsymbol{c}^\calS}\xspace}
\newcommand{\ct}{\ensuremath{\boldsymbol{c}^\mathrm{t}}\xspace}
\newcommand{\co}{\ensuremath{\boldsymbol{c}^\mathrm{o}}\xspace}

\newcommand{\Qs}{\ensuremath{Q^\calS}\xspace}
\newcommand{\Qt}{\ensuremath{Q^\mathrm{t}}\xspace}
\newcommand{\Qo}{\ensuremath{Q^\mathrm{o}}\xspace}

\newcommand{\It}{\ensuremath{\mathbb{I}^\mathrm{t}}\xspace}

\newcommand{\Hh}{\ensuremath{\mathcal{H}}\xspace}
\newcommand{\pVt}{\ensuremath{\partial\Vt}\xspace}
\newcommand{\Itnew}{\mathbb{I}^{n_{\mathrm{t}}}_{\mathrm{new}}}
\newcommand{\so}{\mathfrak{so}(3)}
\newcommand{\se}{\mathfrak{se}(3)}
\newcommand{\mcode}[1]{\lstinline[style=Matlab-editor]!#1!}

\newcommand{\figref}[1]{Fig.~\ref{#1}}

\DeclareMathOperator*{\eig}{eig}

\begin{document}
\begin{frontmatter}

\title{Path Planning for Concentric Tube Robots: a Toolchain with Application to Stereotactic Neurosurgery\thanksref{footnoteinfo}} 

%Constrained Path Planning in Stereotactic Neurosurgery:  Solution Approach via a Homotopy Method

\thanks[footnoteinfo]{This work is funded by the German Research Foundation (DFG, grants FL 989/5-1, OE 253/7-1, SA 773/10-1, WO 2056/11-1).}

\author[First]{Matthias K. Hoffmann} 
\author[First]{Willem Esterhuizen} $^{**}$
\author[Second]{Karl Worthmann}
\author[First]{Kathrin Flaßkamp}

\address[First]{Systems Modeling and Simulation, Saarland University, Germany (e-mail: kathrin.flasskamp/matthias.hoffmann@uni-saarland.de).}
\address[Second]{Institute of Mathematics, Ilmenau University of Technology, Germany (e-mail: willem-daniel.esterhuizen/karl.worthmann@tu-ilmenau.de)}

\begin{abstract}
We present a toolchain for solving path planning problems for concentric tube robots through obstacle fields. First, ellipsoidal sets representing the target area and obstacles are constructed from labelled point clouds. Then, the nonlinear and highly nonconvex optimal control problem is solved by introducing a homotopy on the obstacle positions where at one extreme of the parameter the obstacles are removed from the operating space, and at the other extreme they are located at their intended positions. We present a detailed example (with more than a thousand obstacles) from stereotactic neurosurgery with real-world data obtained from labelled MPRI scans.

\end{abstract}


\end{frontmatter}
%===============================================================================

\section{Introduction}

%The surgery on structures deeply located inside the brain of a patient is a delicate task.\
%While surgery can be performed on the open brain, this has the risk of damaging surrounding areas and leaving the patient with a large portion of the skull removed.\
%The modern, minimally invasive stereotactic neurosurgery instead uses straight cannulae entering the brain from a small hole on the skull.\
%Based on magnetic resonance imaging and computer tomography data, the surgery is arranged using a planning software to reach the target, while avoiding sensitive areas and vessels.\
%Due to distorted anatomy, this task can become difficult or impossible using straight cannulae, so that the surgeon has to rely on open surgery with its disadvantages.\\
%Concentric-tube continuum robots could prove as one solution to the question on how to allow stereotactic neurosurgery in these more difficult scenarios.\ 
%The nesting of bent tubes allows for curved cannulae configurations able to traverse the convoluted brain structures to reach the target area when no straight path is possible.\\

Concentric-tube continuum robots possess great potential to improve stereotactic surgery, due to their ability to trace out curved paths in the body. This is particularly true in neurosurgery where it may be desirable to reach a target area while avoiding various sensitive structures and blood vessels of the brain. The modelling and control of concentric tube robots have received a lot of attention, see  \cite{gilbert2016concentric} for a review. Models of such robots range from those derived from kinematic and/or geometric arguments, such as in \cite{dupont2009design}, \cite{bergeles2015concentric} and \cite{granna2019computer}, to more complicated ones involving physical effects due to bending and torsion of the tubes, as in \cite{webster2006toward} and \cite{rucker2011mechanics}. The papers \cite{greiner2017influence} and \cite{ha2018modeling} aim at also describing certain nonlinear effects.

Much research has been done on path planning for concentric tube robots, with a particular focus on stereotactic surgery. Using a model from \cite{webster2009closed}, the paper by \cite{lyons2009motion} is able to state the path planning problem as a finite dimensional optimization problem. A ``sample-based motion planning'' approach is presented by \cite{torres2011} where the problem is addressed with rapidly exploring roadmaps that use the model by \cite{rucker2011mechanics}. Also using the model from \cite{rucker2011mechanics}, and applying a sample-based approach, the paper \cite{burgner2013computational} tries to maximise the volume reachable by the tube tip subject to a constrained workspace. The paper by \cite{peikert2022automated} addresses the problem by finding paths of connected voxels that the robot can traverse, subject to constraints on its curvature. \cite{flasskamp2019towards} consider obstacle-avoiding path planning in neurosurgery for a tube robot in two dimensions stated as an optimal control problem. The cost functional is a weighted sum of various costs that try to minimise brain damage and error to the target position. \cite{icinco22} also investigate path planning as an optimal control problem with various cost functions, using the model of \cite{rucker2011mechanics}. The paper \cite{leibrandt2017concentric} presents software capable of computing a large set of possible tube configurations as predicted by the kinematic model by \cite{dupont2009design}.

The paper by \cite{sauerteig2022optimal} presents an obstacle-avoiding path planning problem for the model derived in \cite{rucker2011mechanics}. The authors model sensitive brain areas as ellipsoids that need to be avoided and investigate the solutions found when optimising either one of two cost functions (one minimising arc length and one penalising distance to the target set). The numerical experiments were conducted on ``toy data'' to demonstrate the idea.

In the current paper we build on the research done in \cite{sauerteig2022optimal}, using the same model, and present a new toolchain that solves path planning in stereotactic neurosurgery. First, building on the work by \cite{Hackenberg2021}, we present a new approach to fit ellipsoids to labelled point clouds that identify a target area, obstacles and the skull. We then present an approach to solve the path planning problem via a homotopy on the obstacle positions, similar to ideas presented in \cite{bergman2018combining} and \cite{de2022minimum}. We demonstrate the approach with a detailed example using real data labelled by medical professionals. By introducing the homotopy we are able to solve the difficult path planning problem that, due to the presence of over a thousand ellipsoidal obstacles, is initially unsolvable.

The outline of the paper is as follows. Section~\ref{sec:model} summarises the concentric tube robot model and the path planning problem. Section~\ref{sec:constraints} presents the approach of fitting ellipsoids to labelled point clouds and the results for a given MRI data set. Section~\ref{sec:sol} covers the details of the homotopy applied to the computed obstacle positions and the results for the given case-study. Finally, Section~\ref{sec:conclusion} concludes the paper, with recommendations for future research.

\subsection*{Notation and concepts from rigid-body motion}

Given a vector $\x\in\mbbR^n$ and a real symmetric positive definite matrix $Q\in\mathbb{R}^{n\times n}$, $Q\succ 0$, we let $\Vert \x \Vert_{Q}^2 := \x^\top Q \x$ be the $Q$-weighted Euclidean norm. An ellipsoidal set in $\mbbR^3$ with centre $\boldc\in\mbbR^3$ and matrix $Q\in\mbbR^{3\times 3}$, $Q\succ 0$, is denoted by $\Ell(\boldc,Q) := \{\x\in\mbbR^3 : \Vert \x - \boldc \Vert_Q^2 \leq 1\}$. It contains all the points with a $Q$-distance from $\boldc$ less or equal to 1. %The vectors $\0_n\in\mbbR^n$ and $\1_n\in\mbbR^n$ denote the $n$-dimensional column vectors of zeros and ones, respectively. A finite index set running from $1$ to $K\in\mathbb{N}$ is denoted $\mathbb{I}^K := \{1,2,\dots,K\}$. 
The vector $\0_n\in\mbbR^n$ denotes the $n$-dimensional column vector of zeros. A finite index set is denoted $\mathbb{I}^{K} = \{1,2,\dots,K\}$ with $K\in\mathbb{N}$.
%Throughout the paper, we use standard notation for the well-known Lie groups and Lie algebras associated with rigid-body motions, see for example \cite[App. A]{murray2017mathematical}. 
The special orthogonal group on $\mbbR^3$ is denoted by $\SO=\{\M\in\mbbR^{3\times 3}:\M \M^\top = I,\,\,\mathrm{det}(\M) = +1  \}$, and its associated Lie algebra is denoted by $\so=\{\M\in\mbbR^{3\times 3} : \M^\top = -\M\}$. 
Given a vector $\x:=(x_1,x_2,x_3)^\top\in \mbbR^3$ the wedge operator, $^{\wedge}:\mbbR^3 \rightarrow \so$, produces a skew-symmetric matrix,
\[
	\hat \x= \left(\begin{array}{ccc} 0 & -x_1 & x_2 \\ x_3 & 0 & -x_1 \\ -x_2 & x_1 & 0 \end{array}  \right)
\] 
and the vee operator, $^{\vee}:\so\rightarrow \mbbR^3$ denotes the inverse of the wedge. The special Eucliden group on $\mbbR^3$ is denoted by $\SE := \mbbR^3 \times \SO$, which may be identified with the space of all $4\times 4$ matrices of the form $g = \left[\begin{array}{cc} \R & \p \\ \0_3^\top & 1  \end{array}\right]$, $\R\in\SO$, $\p\in\mbbR^3$. Its associated Lie algebra is denoted $\se:= \{(\p, \hat \boldomega): \p\in\mbbR^3, \hat \boldomega\in\so\}$. Similarly, the wedge operator $^{\wedge}:\mbbR^6\rightarrow \se$ maps the \emph{twist coordinates} $\xi = (\boldsymbol{v}^\top, \boldsymbol{\omega}^\top)^\top:=(v_1,v_2,v_3,\omega_1,\omega_2,\omega_3)^\top\in \mbbR^6$ to a matrix (the \emph{twist}), 
\[
	\hat \xi= \left(\begin{array}{cc}  \hat{\boldsymbol{\omega}} & \boldsymbol{v}\\ \boldsymbol{0}_3^\top & 0  \end{array}  \right),
\] 
and $^{\vee}:\se\rightarrow \mbbR^6$ denotes its inverse, see for example \cite[App. A]{murray2017mathematical}.
%Consider a rigid body undergoing a transformation $g(s) = \left[\begin{array}{cc} \R(s) & \p(s) \\ \0_3^\top & 1  \end{array}\right]\in \SE$, $s\in\mbbR$ denoting a parameter. Then, $g^{-1}(s)\dot{g}(s) = 
%\left[
%\begin{array}{cc} \R^\top(s)\dot{\R}(s) & \R^\top(s) \dot{\p}(t) \\ \0_3^\top & 0,
%\end{array}
%\right]\in\se$, where the dot denotes differentiation with respect to the parameter. This last matrix is a twist. Thus, if we specify the twist coordinates $\xi(s) = (\boldsymbol{v}(s)^\top, \boldsymbol{\omega}(s)^\top)^\top$ for all $s$ of interest, we obtain:
%\begin{align*}
%	\R^\top(s)\dot{\R}(s) & = \hat{\boldsymbol{\omega}}(s),\\
%	\R^\top(s)\dot{\p}(s) & = \boldsymbol{v}(s).
%\end{align*}

%Then, the curve traced out by a point on the body relative to an inertial frame, $\q(t)\in\mbbR^3$, for $t\in T\subseteq \mbbR$ satisfies the following differential equation:
%\[
%	\left[
%		\begin{array}{c}
%			\dot{\q}(t)\\
%			0
%		\end{array}
%	\right]
%	= 
%	\hat{\xi}(t) 
%	\left[
%	\begin{array}{c}
%		\q(t)\\
%		1
%	\end{array}
%	\right]
%\]
%where the twist is given by,
%\[
%	\hat{\xi}(t) = g^{-1}(t) \dot{g}(t) =  
%	\left[
%		\begin{array}{cc} \R^\top(t)\dot{\R}(t) & \R^\top(t) \dot{\p}(t) \\ \0_3^\top & 0,
%		\end{array}
%	\right].
%\]
%
%\red
%
%Thus, given twist coordniates $\xi=(v,\omega)$, we have:
%\[
%\dot{g}(t) =  g(t) \hat{\xi}(t),
%\]
%and so $\omega = R T R dot$ and $v = R T p dot$...
%
% \blk

\section{Tube Robot Model \& Path Planning Problem}\label{sec:model}

This section first covers the concentric-tube robot model considered in this paper, as presented in \cite{rucker2010geometrically} and \cite{rucker2011mechanics}, see also the paper by  \cite{sauerteig2022optimal}. Then, the path planning problem is presented.

\subsection{Concentric Tube Robot Model}

\begin{figure}
	\begin{center}
		\includegraphics[width=8.4cm]{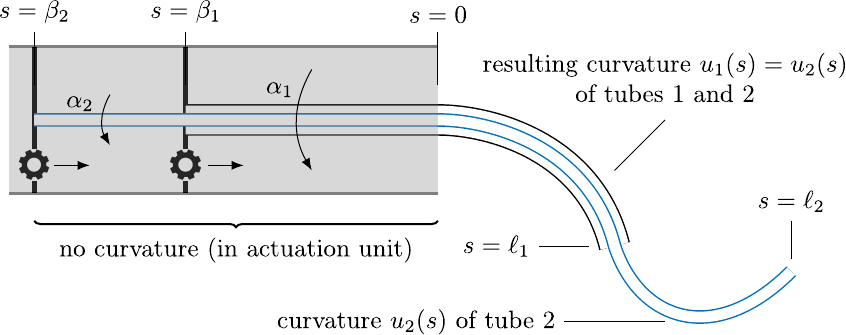}    % The printed column width is 8.4 cm.
		\caption{Example of two concentric tubes with pre-curvature $\boldu_i^\star$ extended by $\beta_i$ and rotated by $\alpha_i$, $i=1,2$ inside the actuation unit. This results in a path $\p(s)$ with curvature $\boldu(s)$, $s\in[\bar{\beta},\bar{\ell}]$ denoting arc length. The tubes are straight for $s\leq 0$, and coincide for $s\in [0,\ell_1]$. Image taken from \cite{sauerteig2022optimal}.} 
		\label{fig:tubes}
	\end{center}
\end{figure}

Consider $n\in\mathbb{N}$ concentric tubes and let $i\in\{1,2,\dots,n\}$ indicate the tube index running from the outer tube to the inner tube (for example, tube 3 runs within tube 2, which runs within tube 1). Each tube has a total length of $L_i\in\mbbR_{\geq 0}$ with a part contained inside the actuation unit, for $s\in[\beta_i,0]$, $\beta_i \leq 0$, and a part that extends outside the actuation unit, for $s\in[0,\ell_i]$, where $s$ is the arc-length, see Figure~\ref{fig:tubes}. Thus, $L_i = \ell_i - \beta_i$.
Each tube by itself (that is, not yet inserted into any other tube) traces out a curve in $\mbbR^3$, denoted $\p_i^{\star}(s)\in \mbbR^3$. Attached to each tube is a right-handed coordinate frame, continuous with respect to arc-length, with the $z$-axis tangent to the curve's velocity vector. Thus, each tube has an associated homogeneous transformation, $g_i^{\star}(s)\in \SE$,
\[
	g_i^{\star}(s) = 
	\left[
		\begin{array}{cc}
			\R_i^{\star}(s) & \p_i^{\star}(s)\\
			\0_3^\top & 1
		\end{array}
	\right],\quad i=1,2,\dots,n,
\]
$s\in [\beta_i,\ell_i]$. Each tube's \emph{pre-curvature} as a function of arc-length is specified by $\boldu _i^{\star}(s) = (\boldu_{ixy}^{\star}(s), u_{iz}^{\star}(s))^\top\in\mbbR^3$, $ \boldu_{ixy}^{\star}(s) = (u_{ix}^{\star}(s), u_{iy}^{\star}(s))^\top\in\mbbR^2$ and satisfies,
\begin{equation*}
	\boldu _i^{\star}(s) = \left((\R_i^{\star}(s))^\top \dot{\R}_i^{\star}(s) \right)^{\vee}\in\mbbR^3.
\end{equation*}
Furthermore, each tube's cross-section is an annulus with constant inner and outer diameter, $\rho_i^{\mathrm{i}}$ and $\rho_i^{\mathrm{o}}$, respectively; $I_i$ denotes its second moment of area (which is constant) about the $x$ or $y$ axis; $E_i$ denotes its Young's modulus; $G_i$ denotes its shear modulus; and $J_i = 2I_i$ denotes its polar moment.

When the tubes are inserted into one another they interact and deform, tracing out a curve with position $\p_i(s)$ and rotation $\R_i(s)$. Let this transformation be denoted,
\[
g_i(s) = 
	\left[
		\begin{array}{cc}
			\R_i(s) & \p_i(s)\\
			\0_3^\top & 1
		\end{array}
	\right],\quad i=1,2,\dots,n.
\]
It is assumed that an outer tube does not extend beyond any of its inner tubes. The effects of various phenomena, such as friction, hysteresis, etc. are ignored, see \cite{rucker2011mechanics} for full assumption details. 

\sloppy
Where the tubes overlap, their positions coincide. However, they are free to rotate about their local $z$-axes. This rotation is denoted by $\psi_i(s)$, which can be shown to satisfy the following differential equations, for each $i$,
\begin{align*}
	\dot{\psi}_i(s) & = u_{iz}(s), \\
	\dot{u}_{iz}(s) & = A_i \sum_{j \in T(s)} E_j I_j (\boldu_{ixy}^{\star}(s))^\top B_{\psi_{ij}}(s)\boldu_{jxy}^{\star}(s),
\end{align*}
for $s\leq \ell_i$. Here, $A_i \coloneqq \frac{E_i I_i}{(EI)(s)G_i J_i }$, $(EI)(s) \coloneqq \sum_{i\in T(s)} E_i I_i$ and $T(s)\coloneqq \{i\in\{1,2,\dots,n\}: s \leq \ell_i\}$ denotes the tube indices of length at least $\ell_i$. Moreover,
\[
	B_{\psi_{ij}}(s) \coloneqq 
	\left[
		\begin{array}{cc}
			\sin(\psi_i(s) - \psi_j(s)) & -\cos(\psi_i(s) - \psi_j(s))\\
			\cos(\psi_i(s) - \psi_j(s)) & \sin(\psi_i(s) - \psi_j(s))
		\end{array}
	\right].
\]
As with the pre-curvatures, we let $\boldu _i(s) = (\boldu_{ixy}(s), u_{iz}(s)) =  (u_{ix}(s), u_{iy}(s), u_{iz}(s))^\top\in\mbbR^3$. Furthermore, the remaining components of the resulting curvature of each tube is given by the following algebraic equation,
\[
	\boldu_{ixy}(s) = \frac{1}{(EI)(s)}\sum_{j\in T(s)} C_{\psi_{ij}}(s)E_j I_j \boldu^{\star}_{jxy}(s),
\]
\fussy where
\[
C_{\psi_{ij}}(s) \coloneqq 
	\left[
		\begin{array}{cc}
			\cos(\psi_j(s) - \psi_j(s)) & -\sin(\psi_j(s) - \psi_i(s))\\
			\sin(\psi_j(s) - \psi_i(s)) & \cos(\psi_j(s) - \psi_i(s))
		\end{array}
	\right].
\]
The displacement and rotation of a tube inside the actuation unit, $\beta_i\in\mbbR_{\leq 0}$ and $\alpha_i\in\mbbR$, respectively, produces a tortion on each tube over the section $[\beta_i, 0]$ (inside the actuation unit) resulting in,
\[
	\psi_i(0) = \alpha_i - \beta_i u_{iz}(0).
\]
Furthermore, it is assumed that there is no external load on the tube, and thus,
\[
	u_{iz}(\ell_i) = 0.
\]
\sloppy Recall the assumption that each tube's local $z$-axis points in the direction of the curve's velocity vector. Thus, also taking the relation $\boldu _i(s) = \left((\R_i(s))^\top \dot{\R}_i(s) \right)^{\vee}$ into account, we identify the twist coordinates, $\xi_i(s) = (\boldsymbol{e}_3^\top, \boldu_i(s)^\top)$, $\boldsymbol{e}_3 \coloneqq (0,0,1)^\top$. The matrix $g^{-1}(s)\dot{g}(s) = 
\left[
\begin{array}{cc} \R^\top(s)\dot{\R}(s) & \R^\top(s) \dot{\p}(t) \\ \0_3^\top & 0,
\end{array}
\right]\in\se$ is a twist. Thus, we obtain, for each $i\in\{1,2\dots,n\}$,
\begin{align*}
	\dot{\R}_i(s) & = \R_i(s)\hat{\boldu}_i(s),\\
	\dot{\p}_i(s) & = \R_i(s)\boldsymbol{e}_3,
\end{align*}
for all $s\in[\beta_i,\ell_i]$.

\subsection{The Path Planning Problem}

Recall that all the tubes' positions coincide and that no tube extends beyond any of its inner tubes. Thus, it suffices to only consider the evolution of $\p_n$. Our goal is to choose, for $i = 1,2,\dots,n$, the precurvatures, $\boldu_i^{\star}$ (which we assume constant for all $s\in[\beta_i,\ell_i]$); the tube lengths, $L_i$; the tube inner and outer diameters, $\rho_i^{\mathrm{i}}$ and $\rho_i^{\mathrm{o}}$; the actuator parameters, $(\alpha_i, \beta_i)$; the initial position of the inner-most tube on the skull, $\p_n^0$; and its initial orientation $\R_n^0$, such that the resulting curve traced out by the inner tube reaches a target area while avoiding a number of obstacles enclosing sensitive brain areas. The tube is also not allowed to cross from one brain hemisphere to the other.

Let $\boldu^\star \coloneqq (\boldu^{\star \top}_1, \boldu^{\star \top}_2,\dots,\boldu^{\star \top}_n)^\top\in\mbbR^{3n}$, $(\boldsymbol{\rho}^{\mathrm{i}}, \boldsymbol{\rho}^{\mathrm{o}})  := ({\rho}^{\mathrm{i}}_1,{\rho}^{\mathrm{i}}_2,\dots,{\rho}^{\mathrm{i}}_n, {\rho}^{\mathrm{o}}_1,{\rho}^{\mathrm{o}}_2,\dots,{\rho}^{\mathrm{o}}_n)^\top\in\mbbR^{2n}$, $(\boldalpha^\top, \boldbeta^\top)^\top \coloneqq (\alpha_1, \alpha_2, \dots, \alpha_n, \beta_1,\beta_2,\dots,\beta_n)^\top\in \mbbR^{2n}$, $\boldL \coloneqq (L_1,L_2,\dots, L_n)^\top\in\mbbR^n$, and let the decision space be denoted by,
\[
	\mathbb{D} \coloneqq \{ (\boldu^\star, \boldL, \boldsymbol{\rho}^{\mathrm{i}}, \boldsymbol{\rho}^{\mathrm{o}}, \boldalpha, \boldbeta, \p_n^0, \R_n^0)\in\mbbR^{8n+ 3}\times\SO\}.
\]
Our path planning problem (PPP) may be expressed as follows,
\begin{align}
	(\PPP)&\quad\quad\quad\quad\min_{\boldd\in\mathbb{D}} \quad  J(\boldd)\nonumber\\
	\mathrm{s.t.} \quad &\mathrm{For\ all\ }s\in[0,\ell_n]\mathrm{\ and\ all\ }i\in\{1,2,\dots,n\}: \nonumber \\
	\quad & \p_n(s) \in (\mathcal{S} \cap \Hh)\setminus\mathcal{E}. \nonumber\\
	\quad &\dot{\p}_n(s) = \R_n(s)\bolde_3, \label{PPP_eq_1} \\
	\quad & \dot{\R}_n(s) = \R_n(s)\hat{\boldu}_n(s), \\
	& \dot{\psi}_i(s) = u_{iz}(s),  \\
	& \dot{u}_{iz}(s) = A_i \sum_{j \in T(s)}  E_j I_j (\boldu_{ixy}^{\star})^\top B_{\psi_{ij}}(s) \boldu_{jxy}^{\star}, \\
	& \boldu_{ixy}(s) = \frac{1}{(EI)(s)}\sum_{j\in T(s)} C_{\psi_{ij}}(s)E_j I_j \boldu^{\star}_{jxy}.\\
	\quad & \p_n(0) = \p_n^0 \in \partial \mathcal{S}, \\
	\quad & \R_n(0) = \R_n^0 \in \SO, \\
	\quad & \p_n(\ell_n) \in \Pf,\\
	\quad & 0 < L_1 \leq L_2 \leq \dots \leq L_n.\\
	\quad & u_{iz}(\ell_i) = 0,\\
	\quad & \psi_i(0) = \alpha_i - \beta_i u_{iz}(0),\\
	\quad & \ell_i = L_i + \beta_i,\\
	\quad & 0 \leq \alpha_i < 2\pi,\\
	\quad & -L_i \leq \beta_i \leq 0. \label{PPP_eq_2}
\end{align}

\fussy
The cost function is chosen to minimise arc length, thus, $J: \mathbb{D} \rightarrow \mbbR_{\geq 0}$ is given by,
\[
J(\boldd) \coloneqq \ell_n.
\]
The skull and brain are modelled by an ellipsoid, $\calS=\Ell(\cs,\Qs)$; $\Pf\subset\mbbR^3$ indicates the target set; $\Hh$ is a half space containing the brain hemisphere wherein the target set lies; and $\mathcal{E} \coloneqq \bigcup_{j\in\mathbb{I}^K}\mathcal{E}_j$, with $\mathcal{E}_j = \Ell(\mathbf{c}^{\mathcal{E}_j}, Q^{\mathcal{E}_j})$, and $K\in\mathbb{N}$, denotes the ellipsoidal obstacles. If a solution to the problem exists, we denote it by $\bar{\boldd}\in \mathbb{D}$.

% \red $\R^0_n$ can be constructed in multiple ways, e.g. with Euler angles or quaternions. \blk
\section{Obtaining the Constraints from Labelled Data}\label{sec:constraints}

This section describes how we construct the various constraints appearing in (PPP) from labelled data, which we assume to be point clouds (voxel centres) in $\mbbR^3$ that indicate relevant areas of a patient's brain. These are points that define the skull, $\Vs = \left\lbrace \vs_i \right\rbrace_{i\in\mathbb{I}^{\ns}}$; the target area, $\Vt = \left\lbrace \vt_i \right\rbrace_{i\in\mathbb{I}^{\nt}}$; $K\in\mathbb{N}$ obstacles, $\Vo_j = \left\lbrace \vo_{ij} \right\rbrace_{i\in\mathbb{I}^{\noj}}$, $j=1,2,\dots,K$; and points on a plane that divide the brain into its hemispheres, $\Vh = \left\lbrace \vh_i \right\rbrace_{i\in\mathbb{I}^{\nh}}$.

\subsection{Fitting the hyperplane}

The plane dividing the hemispheres, which we label $\Pi:=\{\x\in\mbbR^3: \x^\top \boldsymbol{h} = 1  \}$, $\boldsymbol{h}\in\mbbR^3$, is found by solving for $\boldsymbol{h}$ in the system of equations,
\[
(\vh_i)^\top \boldsymbol{h} = 1,\quad i\in\mathbb{I}^{\nh},
\]
via least-squares linear regression. If all the target points $\Vt$ are contained in one half space defined by this plane, then we take $\Hh$ to be this half space. Otherwise, if target points appear in both hemispheres, we solve (PPP) twice: once with $\Hh = \{\x: \x^\top \boldsymbol{h} \leq 0\}$, and once with $\Hh = \{\x: \x^\top \boldsymbol{h} \geq 0\}$ and choose the best solution of the two problems.

\subsection{Fitting the skull}

The ellipsoid containing the brain and the skull, $\mathcal{S}=\Ell(\cs,\Qs)$, is found as a best-fit ellipsoid for all the points in $\Vs \cap \Hh$ by solving,
\begin{align*}
	\min_{\cs, \Qs} \quad& \sum_{i \in \mathbb{I}^{\ns}_{\mathrm{new}} } \left( \Vert \vs_i - \cs \Vert_{\Qs}^2 - 1\right)^2\\
	\text{s.t.} \quad & \Qs \succ 0,\\
	& \cs\in\mbbR^3,
\end{align*}
where $\mathbb{I}^{\ns}_{\mathrm{new}} := \left\{i\in \mathbb{I}^{\ns} : \vs_i \in (\Vs \cap \Hh) \right\}$.

By using fewer points in the fitting as we have done here we obtain faster execution times without loss of accuracy. Positive definiteness for a matrix $Q$ can be achieved in the following two ways.\ 
One is the Cholesky decomposition, $Q = GG^\top$, with $G \in \mathbb{R}^{3\times 3}$ a lower triangular matrix and the main diagonal entries greater than 0.\ 
The second, and in our experiments more stable one, is with rotation matrices, $Q = R^\top M R$, with $M \in \mathbb{R}^{3\times 3}$, a diagonal matrix with positive entries on the main diagonal, and $R \in \SO$, described by quaternions.\
This formulation allows for setting constraints of the eigenvalues of $Q$ easily.\ 

\subsection{Fitting the target}

The target area is modelled as an ellipsoid enclosed by the target points, so that it lies completely inside the point cloud of the target.\ 
The target's most exterior points, forming the boundary \pVt of the point cloud can be for example found with MATLAB's \lstinline[style=Matlab-editor]|boundary|-function.\ 
A tight boundary (meaning not just the convex hull) is preferred, so that no part of the ellipsoid reaches outside of the point cloud.\
Let $\delta_\mathrm{MRI}$ denote the minimal distance between any two points of \Vt.
%The target points have a \red minimal distance of $\delta_\mathrm{MRI}$\
If one of the ellipsoid's semi-axes was shorter than $\delta_\mathrm{MRI}$, it could lie between all the target points and thus exit the point cloud, so the eigenvalues of \Qt are constrained to prevent this, resulting in the optimization problem,
\begin{align*}
	\min_{\ct, \Qt} \quad& \sum_{i \in \Itnew} \Vert \vt_i - \ct \Vert_{\Qt}^2 \\
	\text{s.t.} \quad& \eig(\Qt) > \left(\frac{2}{\delta_\mathrm{MRI}}\right)^2\\
	\quad& \Vert \vt_i - \ct \Vert_{\Qt}^2 \ge 1 \quad \forall \ i \in \Itnew,
\end{align*}
where $\Itnew:= \{i\in\It : \vt_i \in\partial \Vt   \}$, so that,
\[
\Pf = \left\{\boldsymbol{x} : \Vert \boldsymbol{x} - \ct \Vert_{\Qt}^2 \le 1\right\}.
\]
The minimisation of the target points' $\Qt$-distance is equivalent to maximising the ellipsoid's volume, the product of the eigenvalues of $M$, but with faster convergence and more consistent results.\

\subsection{Fitting the obstacles}

The same trick from the previous subsection is used in the formulation of enclosing ellipsoids for the obstacles.\ 
For an obstacle set $\Vo_{j}$ we solve,
\begin{align}
	\label{eq:ell_obs}
	\min_{\co_j, \Qo_j} \quad& -\sum_{i \in \mathbb{I}^{\noj}_{\mathrm{new}}:} \Vert \vo_{ij} - \co_j \Vert_{\Qo_j}^2 \nonumber\\
	\text{s.t.} \quad& \Qo_j \succ 0\nonumber\\
	\quad& \Vert \vo_{ij} - \co_j \Vert_{\Qo_j}^2 \le 1 \quad \forall \ i \in \mathbb{I}^{\noj}_{\mathrm{new}}
\end{align}
where $\mathbb{I}^{\noj}_{\mathrm{new}}:= \{i\in \mathbb{I}^{\noj} : \vo_{ij}\in  \Vo_j \cap \Hh\}$. This provides the enclosing ellipsoids,
$\mathcal{E}_j = \left\{ \boldsymbol{x} : \Vert \boldsymbol{x} - \co_j \Vert_{\Qo_j}^2 \ge 1 \right\}$.

\subsection{Division of obstacle points}

When building the obstacle ellipsoids, we are interested in fully covering all of the obstacle points, given by the construction of enclosing ellipsoids, but also not include large volumes of non-obstacle regions.\ 
One major contributor to this problem is if non-connected areas are captured by the same ellipsoid, as it might happen with the k-means algorithm, shown by \cite{Hackenberg2021}.\ 
Applying DBSCAN, a density-based clustering algorithm by \cite{Ester1996dbscan}, beforehand can alleviate this flaw.\
The distance $\delta_\mathrm{MRI}$ is used as the search radius for DBSCAN.
To evaluate how well the ellipsoids cover the obstacle points and the surrounding area, the coverage is calculated as the number of obstacle points divided by the number of grid points on the MRI grid inside the cluster.\
Algorithm~\ref{algo:clustering} shows the pseudo-code describing the steps.

\begin{algorithm}
	\label{algo:clustering}
	\caption{Clustering of obstacle points}
	\begin{algorithmic}[1]
		\Require{obstacle points \vo and threshold $c_\mathrm{th}$}
		\State{Cluster \vo using DBSCAN to get $\vo_i$.}
		\For{each $i$}
			\While{true}
				\State{Cluster $\vo_i$ with k-means to get $\vo_{ij}$.}
				\For{each $j$}
					\State{Solve \eqref{eq:ell_obs} with $\vo_{ij}$.}
				\EndFor
				\If{$\frac{\text{\# MRI grid points in ellipsoids}}{\text{\# obstacle points in ellipsoid}} \le c_\mathrm{th}$}
					\State{Increase the number of clusters for k-means.}
				\Else
					\State{Exit the while loop.}
				\EndIf
			\EndWhile
			\State{Add the $j$ ellipsoids to the obstacle set.}
		\EndFor
	\end{algorithmic}
\end{algorithm}

%From an initial point the algorithm searches for neighbours closer than $\delta_\mathrm{DBSCAN}$, adds them to the cluster and repeats with their neighbours until no further neighbour can be located.\ 
%In that case, a new cluster is started until all points are matched to a cluster or even are marked as outliers.\ 
%As the obstacle points follow the MRI grid, we have a great guess for $\delta_\mathrm{DBSCAN} = \delta_\mathrm{MRI}$, a valuable information for DBSCAN.\ 
%Each of the resulting DBSCAN clusters is clustered with k-means.\
%To evaluate how exactly the ellipsoids cover a DBSCAN cluster we use the ratio of grid points covered by all the corresponding ellipsoids to points in the cluster.\
%If this ratio is below a user-defined threshold the number of k-means clusters is increased.\ 
%The $\Vo_j$ are formed from the points contained in the k-means clusters.
%As we approximate the point clouds with ellipsoids, we expect to cover volumes larger than the voxel volume.

\subsection{Numerical results}
\begin{figure*}[htb]
	\begin{minipage}{0.28\textwidth}
		\centering
		\includegraphics[height=51mm]{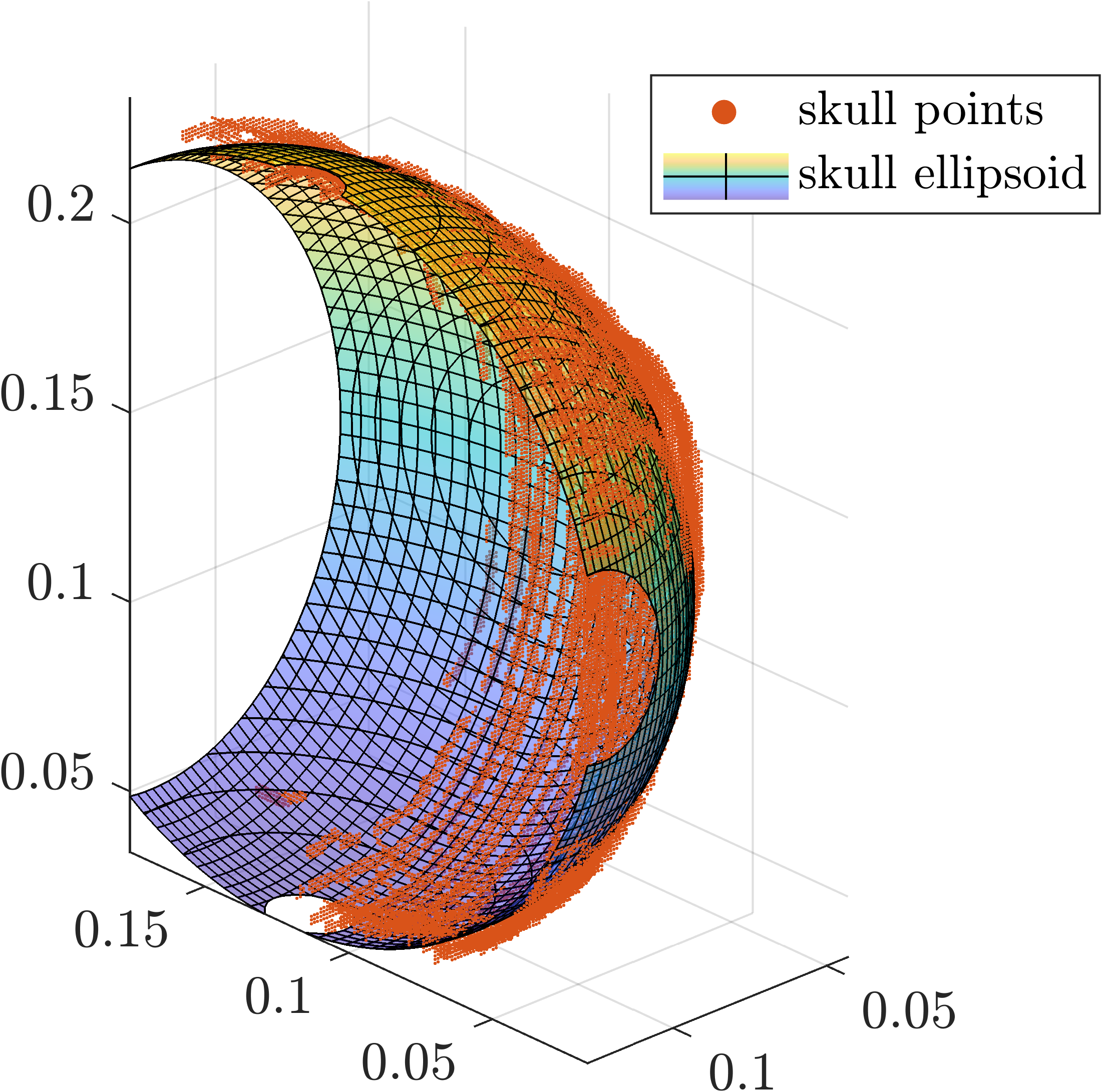}
		\caption{A best-fit ellipsoid is used to model the skull around the allowed hemisphere.\ An ellipsoid seems to be a good approximation for the skull. \vspace{1\baselineskip} }
		\label{fig:skull_opt}
	\end{minipage}
	\vspace{0.03\textwidth}
	\begin{minipage}{0.3\textwidth}
		\centering
		\includegraphics[height=51mm]{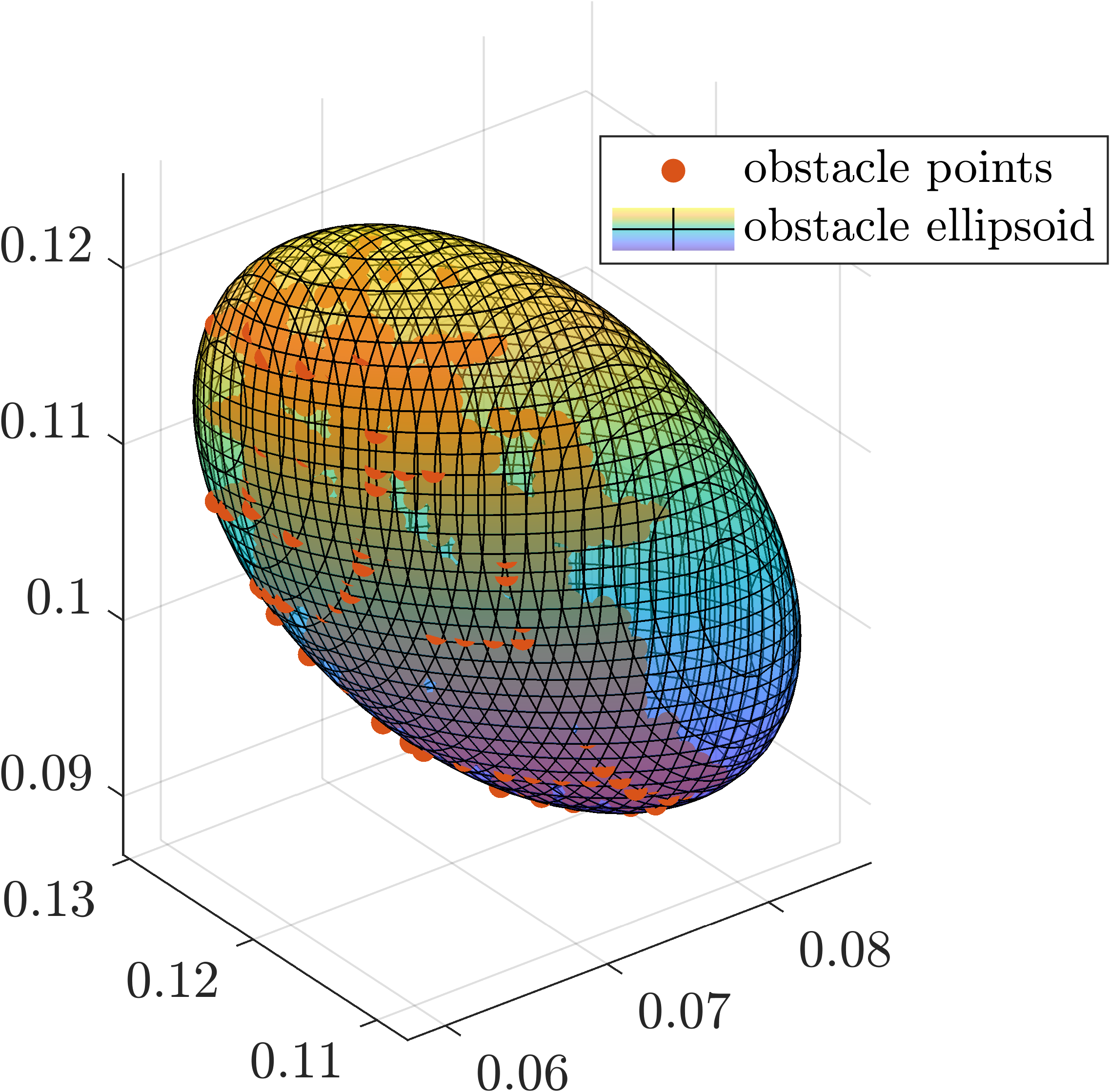}
		\caption{Enclosing ellipsoid as an approximation for the obstacles.\ Parts of the ellipsoid contain no obstacle points, where the ellipsoid is disallowing non-critical parts of the brain.}
		\label{fig:obs_opt}
	\end{minipage}
	\vspace{0.03\textwidth}
	\begin{minipage}{0.36\textwidth}
		\centering
		\includegraphics[height=52mm]{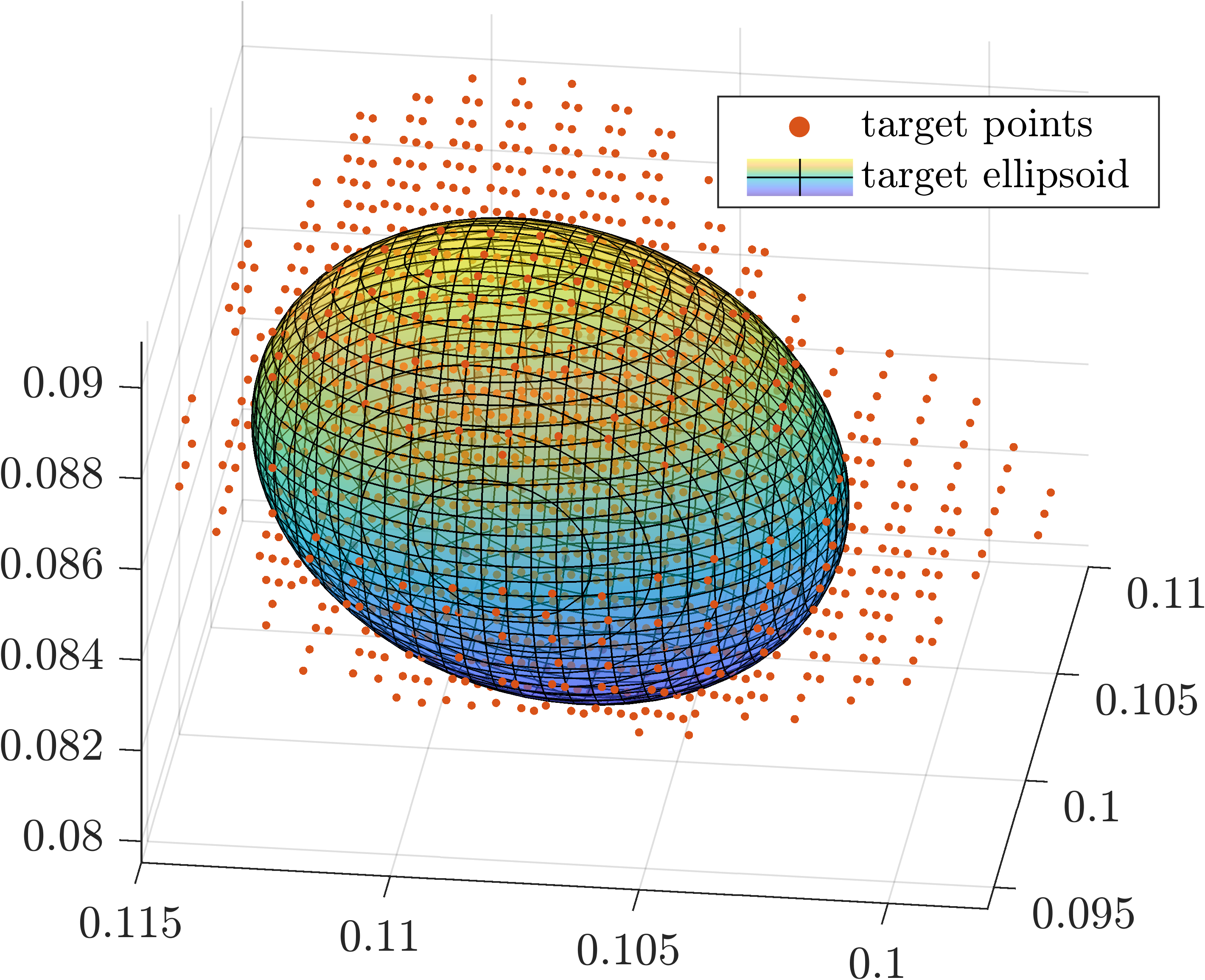}
		\caption{The target volume is modelled as an ellipsoid enclosed by the point cloud.\ The ellispoid is an underestimation to ensure that the cannula enters the target. \vspace{1\baselineskip} }
		\label{fig:target_opt}
	\end{minipage}
\vspace{-1.2cm}
\end{figure*}

Figures~\ref{fig:skull_opt}-\ref{fig:target_opt} show the optimisation results for the three different types of ellipsoids.\ 
In \figref{fig:skull_opt}, only the skull points in the admissible half-space are displayed.\ 
It shows that an ellipsoid is a good approximation for the skull.\ 
As already mentioned before, the ellipsoid in \figref{fig:obs_opt} covers a larger volume than the voxels marked as obstacles.\
This can be seen as there are parts in the ellipsoid without obstacle points.\ 
The target modelling is successful as well, still the ellipsoid underestimates the size of the target region.\ 
Overall, this problem has 1179 obstacle ellipsoids resulting in a convoluted obstacle field, partly shown in \figref{fig:opt_solved_ellipsoids}.\ 

\section{Solution Approach via Homotopy}\label{sec:sol}

The path planning problem we consider is nonlinear and highly nonconvex. With realistic data (where hundreds of ellipsoids may be present) numerical solvers often do not converge and even if they do, need very long execution times. Thus, using similar ideas from \cite{bergman2018combining} and \cite{de2022minimum}, we introduce a homotopy on some of the obstacles' positions, where at the one extreme they are removed from the skull's interior, and at the other they are located at their original position. We then iteratively solve relaxed problems, using the solution with a current parameter as the initial guess for the solution with the next parameter (see Algorithm~\ref{eq:alg_1}).

For a given \emph{homotopy parameter} $\lambda\in[0,1]$ the relaxed problem reads,
\begin{align*}
	(\PPP^{\lambda})&\quad\quad\quad\quad\min_{\boldd\in\mathbb{D}} \quad  J(\boldd)\nonumber\\
	\mathrm{s.t.} \quad &\mathrm{For\ all\ }s\in[0,\ell_n]\mathrm{\ and\ all\ }i\in\{1,2,\dots,n\}: \nonumber \\
	\quad & \p_n(s) \in (\mathcal{S} \cap \Hh)\setminus\mathcal{E}^{\mathrm{fix}}, \nonumber\\
	\quad & \p_n(s) \in (\mathcal{S} \cap \Hh)\setminus\mathcal{E}^{\lambda}, \nonumber\\
	\quad & \eqref{PPP_eq_1} - \eqref{PPP_eq_2}.
\end{align*}
Here we let
\[
\mathcal{E}^{\mathrm{fix}}= \bigcup_{j\in \mathbb{I}^{\mathrm{fix}}}\calE_j,
\] 
where $\mathbb{I}^{\mathrm{fix}}\subseteq\mathbb{I}^K$ is the index set of ellipsoids we do \emph{not} want to perturb with $\lambda$; and
\[
\mathcal{E}^{\lambda} = \bigcup_{j\in\mathbb{I}\setminus\mathbb{I}^{\mathrm{fix}}}\mathcal{E}_{j}^{\lambda},
\]
where
\[
\mathcal{E}_{j}^{\lambda} = \mathrm{Ell}([1- \lambda]\boldc_j^\mathrm{init} + \lambda \boldc_j^{\mathrm{o}}), \Qo_j),
\]
indicates the ellipsoids that we want to perturb. Again, if a solutions exists, we denote it by $\bar{\boldd}^{\lambda}$. Thus, $\lambda$ moves the chosen ellipsoids along the line segment that connects their original positions, $\boldc_j^{\mathrm{o}}$, with arbitrary user-specified positions, $\boldc_j^\mathrm{init}$, where the entire ellipsoid $\calE_j$ is located outside the skull. The \emph{homotopy algorithm} we employ is presented in Algorithm~\ref{eq:alg_1}.
\begin{algorithm}
	\caption{Homotopy Algorithm}\label{eq:alg_1}
	\begin{algorithmic}[1]
		\Require initial guess $\boldd^0$; step size $\Delta > 0$
		\State $\lambda \gets 0$
		\State $\boldd\gets \boldd^0$
		\While{$\lambda < 1$}
		\State try solve $\PPP^\lambda$ with $\boldd$ as initial guess.
		\If{solution to $\PPP^\lambda$ is found}
		\State $\boldd \gets \bar{\boldd}^\lambda$
		\State $\lambda \gets \min\{1,\lambda + \Delta\}$
		\Else
		\State break while loop
		\EndIf
		\EndWhile
		\State \Return $\boldd$, $\lambda$
	\end{algorithmic}
\end{algorithm}

\subsection{Heuristic for calculating user-specified positions $\boldc_j^\mathrm{init}$}
\begin{figure}[htb]
	\centering
	\fontsize{9}{0}\selectfont
	\def\svgwidth{0.49\textwidth}
	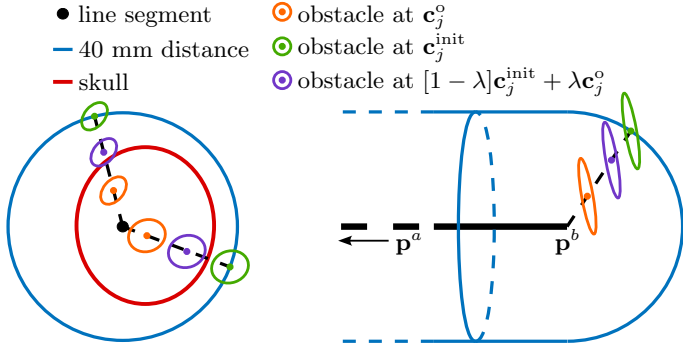
	\caption{The obstacle centres are pushed away perpendicular to the line segment connecting the entry and terminal point of the initial guess.}
	\label{fig:homotopy_shift}
\end{figure}

The problem of finding a feasible cannula configuration is simpler if there are just a few obstacles.\ 
With an increasing number, the problem becomes so hard that the solver does not converge.\ 
As our aim is to shift the obstacles out of the way, one idea is to take an initial guess with very few obstacles for $\mathcal{E}^{\mathrm{fix}}$ (so few that the problem is easily solvable) and solve the problem to get an initial tube configuration. We then let $\p^a = \p(0)$ and $\p^b = \p(\ell_n)$.\
%The 10 ellipsoids closest to the target centre make up $\mathcal{E}^{\mathrm{fix}}$ so that the initial guess will not be arbitrarily bad.\
All the non-fixed obstacles are shifted orthogonally from the line segment connecting $\p^a$ and $\p^b$.\
The factor
\begin{align*}
	a_i = \frac{(\co_i - \p^a)^\top (\p^b - \p^a)}{\Vert(\p^b - \p^a)\Vert_2^2}
\end{align*}
is needed to find the point closest on the line segment
\begin{align*}
	\boldsymbol{l}_i = \left\{\begin{matrix}
		\p^a, &\text{if } a_i < 0\\
		\p^a + a_i(\p^b-\p^a), &\quad \ \ \, \text{if } 0 \le a_i \le 1\\
		\p^b, &\text{if } a_i > 1
	\end{matrix}\right..
\end{align*} 
Thus, the user-specified positions
\begin{align*}
	\boldc_i^\mathrm{init} = \co_i + 0.4\frac{\co_i - \boldsymbol{l}_i}{\Vert \co_i - \boldsymbol{l}_i \Vert_2}, \ \forall \ i \in (\mathbb{I}\setminus\mathbb{I}^{\mathrm{fix}})
\end{align*}
are shifted to a distance of \SI{400}{\milli\meter} from the line segment, which is typically outside of the skull.\
\figref{fig:homotopy_shift} shows the approach on an example.\ 

\subsection{Numerical results}
Next, we display the numerical results of using the homotopy algorithm on a case study.\
CasADi by \cite{Andersson2019} is used to formulate the optimisation problem, while $(\PPP^{\lambda})$ is solved with IPOPT by \cite{Waechter2006}.\ 
Even though interior-point methods are known to be hard to warm-start (see for example \cite{John2008ipwarmstart}), we show that they work well in this problem. The following table contains the values of the IPOPT parameters we changed.

\begin{tabular}{l|c}
	\mcode{warm_start_init_point} & \mcode{yes}\\
	\mcode{mu_init} & \mcode{1e-8}\\
	\mcode{warm_start_mult_bound_push} & \mcode{1e-10}\\
	\mcode{warm_start_slack_bound_push} & \mcode{1e-10}\\
	\mcode{warm_start_bound_push} & \mcode{1e-8}\\
	\mcode{warm_start_bound_frac} & \mcode{1e-8}\\
	\mcode{warm_start_slack_bound_frac} & \mcode{1e-10}\\
\end{tabular}

The given case-study is intended to be solved using three tubes, with $u_{y,i}^\star = 0$, $L_i = \SI{600}{\milli\meter} \ \forall \ i\in[1,3]$, and a maximal outer diameter of the outer-most tube of $\rho^\mathrm{o}_1 = \SI{8}{\milli\meter}$.
The length of the inner tube is minimised over the tube actuations $\boldalpha$ and $\boldbeta$; tube lengths, $\mathbf{L}$; the tube curvatures $u_{i,x}^\star$; the inner and outer tube diameters, $\rho_i^{\mathrm{i}}$ and $\rho_i^{\mathrm{o}}$; the initial condition $\p_n^0$ and the initial rotaion,
$$\R^0_n = \begin{bmatrix}
	1 - 2(q_2^2-q_3^2) & 2(q_1 q_2 - q_0 q_3) & 2(q_1 q_3 + q_0 q_2)\\
	2(q_1 q_2 + q_0 q_3) & 1 - 2(q_1^2-q_3^2)  & 2(q_2 q_3 - q_0 q_1)\\
	2(q_1 q_3 - q_0 q_2) & 2(q_2 q_3 + q_0 q_1) & 1 - 2(q_1^2-q_2^2)
\end{bmatrix},
$$
where $\q^0_n = [q_0, q_1, q_2, q_3]$ is a quaternion.

The original problem $(\PPP)$ with its 1179 obstacles is not solvable out-of-the-box as the solver will not converge. However, a solution is found with the homotopy algorithm with $\Delta = 0.1$. The following tables show the resulting values for the tube parameters and decision variables of tubes 2 and 3, as the algorithm finds that the outer most tube is not necessary, resulting in $\ell_1=0$.
\begin{flushleft}
	\begin{tabular}{l|c|c|c|c}
		Variable& $\ell_2$ & $\ell_3$ & $u_{x,2}^\star$ & $u_{x,3}^\star$ \\
		Value & \SI{13.2}{\milli\meter} & \SI{40.1}{\milli\meter} & \SI{9.99}{\per\meter} & \SI{4.55}{\per\meter}\\
		\hline
		Variable\Tstrut & $\alpha_2$ & $\alpha_3$ & $\beta_2$ & $\beta_3$\\
		Value &1.5405 & 4.7425 & \SI{-0.587}{\meter} & \SI{-0.56}{\meter}\\
		\hline
		Variable\Tstrut & $\rho^\mathrm{i}_1$ & $\rho^\mathrm{o}_1$ & $\rho^\mathrm{i}_2$ & $\rho^\mathrm{o}_2$ \\
		Value & 0.004 & 0.0053 & 0.0013 & 0.0026\\
		\hline
		Variable\Tstrut& \multicolumn{4}{c}{$\p^0_n$} \\ Value & \multicolumn{4}{c}{$\left[0.0967, 0.0769, 0.0536\right]^\top$} \\ \hline
		Variable\Tstrut& \multicolumn{4}{c}{$\q_n^0$} \\ Value & \multicolumn{4}{c}{$\left[0.6983, 0.1727, -0.3151, -0.6190\right]^\top$}
	\end{tabular}
\end{flushleft}

\begin{figure}[htb]
	\centering
	\includegraphics[scale=0.6]{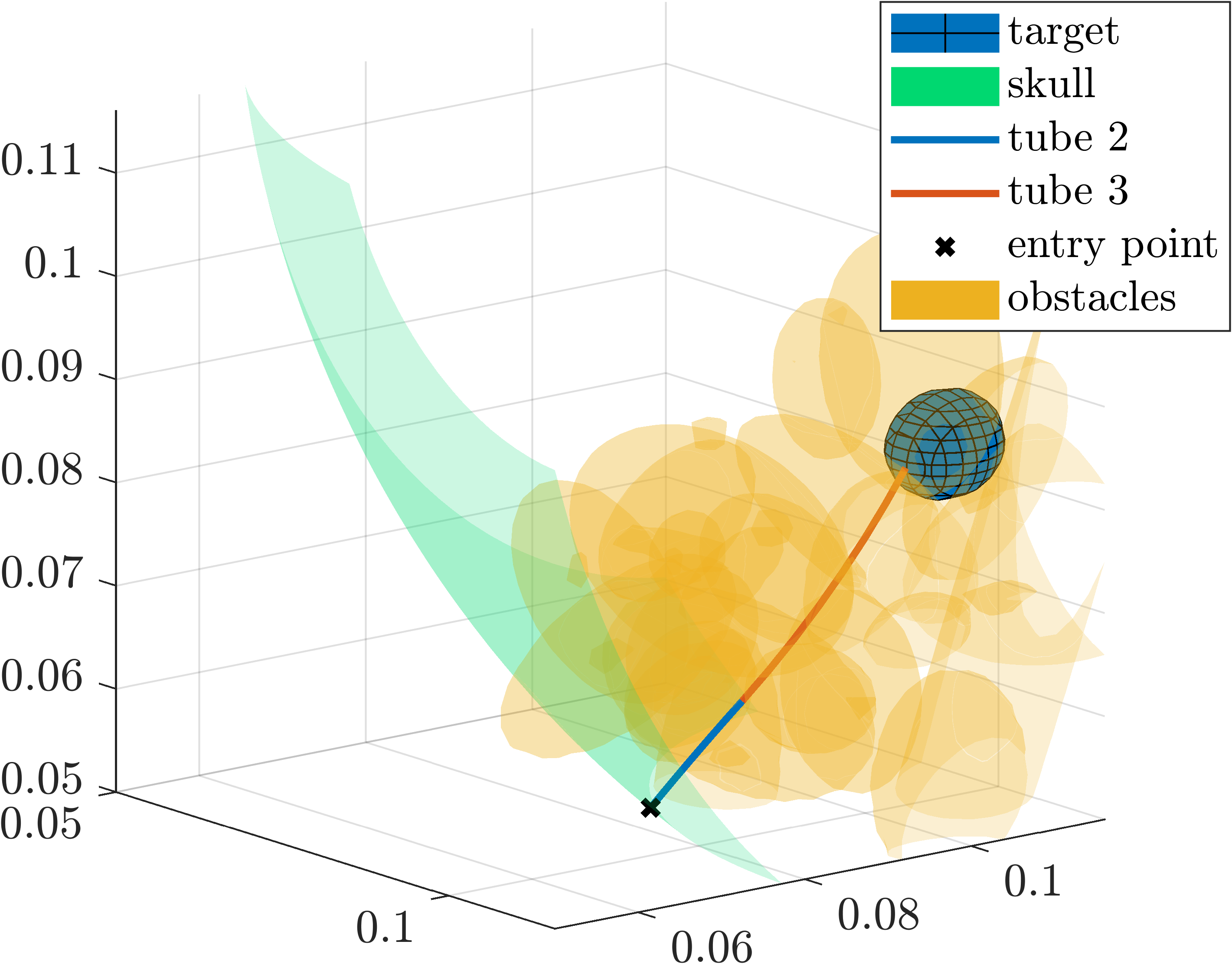}
	\caption{The solution to the given case-study reaching from the skull to the target position.\ The 45 ellipsoids closest to the cannula are displayed.}
	\label{fig:opt_solved_ellipsoids}
\end{figure}

\figref{fig:opt_solved_ellipsoids} shows the optimised cannula navigating through the field of the 45 closest obstacles, starting on the skull and entering the target area.

\begin{figure}[htb]
	\centering
	\includegraphics[scale=0.7]{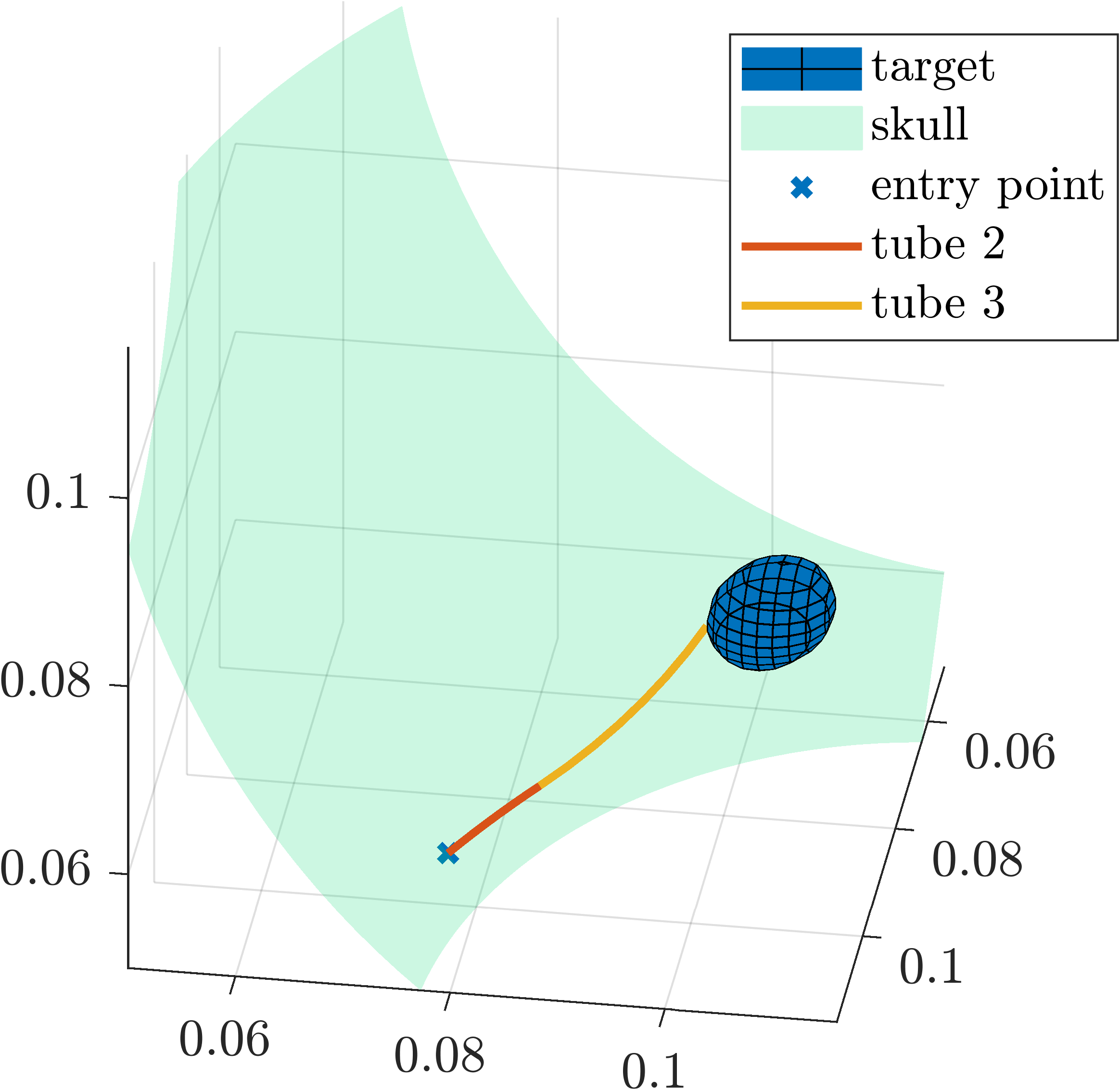}
	\caption{The solution without the obstacles.}
	\label{fig:opt_solved}
\end{figure}

The same solution is displayed in \figref{fig:opt_solved} without the obstacles to show the resulting curvature of the cannula.\

\section{Conclusion}\label{sec:conclusion}

In this work, we tackle the challenging task of optimally planning concentric tube continuum robot dimensioning and actuation in geometrically constrained spaces and present our findings in the real-world example of stereotactic neurosurgery.\ 
From labelled MRI or CT scans as a starting point, we formulate easily solvable optimisation problems for condensing the point clouds to a set of elliptical constraints.\ 
For the large number of points belonging to the obstacles, we propose an algorithm combining DBSCAN and the k-means-algorithm for separating and clustering the obstacles.\ 
By using homotopy methods, we enable the solution of the path planning problem, unsolvable out-of-the-box.\ 
The success of these methods is underlined by the application to the real-world case-study and the corresponding numerical results.

There are many areas on which future research may focus. First, the homotopy algorithm is not guaranteed to converge as $\lambda$ tends to 1 and conditions on the problem data that imply this should be investigated. Second, in moving the obstacles with the homotopy parameter we did not take care to prevent topological changes in the free space where the tubes may manoeuvre. Though we were able to find a solution to the problem we considered, these topological changes may result in discontinuities of the optimal paths from one iteration of $\lambda$ to the next, which could cause infeasibility. Third, by using the homotopy approach the final path may in fact only be locally optimal: a result of ``slightly'' perturbing previously-found solutions. More involved perturbations can result in global optima.

%
%\begin{ack}
%\red Place acknowledgments here. \blk
%\end{ack}

\bibliography{ifacconf}             % bib file to produce the bibliography
                                                     % with bibtex (preferred)
    
\end{document}